\documentclass[11pt,twocolumn,abstract]{scrartcl}

\usepackage[utf8]{inputenc}
\usepackage[T1]{fontenc}

\usepackage{color}
\usepackage{xcolor}

\usepackage{pgfplots}
\pgfplotsset{compat=1.18, compat/show suggested version=false}
\usepackage{pgfplotstable}
\usepgfplotslibrary{colorbrewer}

\usepackage{mathtools}
\usepackage{amsmath,amsfonts,amssymb}
\allowdisplaybreaks[4]

\usepackage{bm}
\usepackage{siunitx}
\usepackage[ISO]{diffcoeff}
\usepackage{cite}

\usepackage{booktabs}

\pgfplotscreateplotcyclelist{scatter markers}{
    every mark/.append style={solid,fill=\pgfplotsmarklistfill, scale=0.80},mark=square*\\
    every mark/.append style={solid,fill=\pgfplotsmarklistfill, scale=1},mark=triangle*\\
    every mark/.append style={solid,fill=\pgfplotsmarklistfill, scale=0.81},mark=*\\
}

\pgfplotscreateplotcyclelist{scatter markers alt}{
    every mark/.append style={solid,fill=\pgfplotsmarklistfill, scale=0.80},mark=square*\\
    every mark/.append style={solid,fill=\pgfplotsmarklistfill, scale=1},mark=triangle*\\
    every mark/.append style={solid,fill=\pgfplotsmarklistfill, scale=0.85}, mark=pentagon*\\
    every mark/.append style={solid,fill=\pgfplotsmarklistfill, scale=0.81},mark=*\\
}

\pgfdeclareplotmark{arrow}{ \draw[line width=0.6pt,<-] (0,-0.7mm) -- ++(0,1.7mm); }

\pgfplotsset{
    stdmarks bw/.style={
        cycle list/Greys-5,
        cycle multiindex* list={
            mark list*\nextlist
            Greys-5\nextlist
            draw=black
        },
        line join=round,
    },
    errormarks bw/.style={
        cycle list/Greys,
        cycle multiindex* list={
            mark list*\nextlist
            Greys-F, Greys-G, Greys-H, Greys-I, Greys-J\nextlist
            {},{},draw=black\nextlist
            linestyles\nextlist
        },
        mark options={solid,draw=black,scale=0.7,line width=0.2pt},
        line join=round,
    },
}

\pgfplotsset{
    stdmarks/.style={
        cycle list/Set2,
        cycle multiindex* list={
            mark list*\nextlist
            Set2\nextlist
            draw=black
        },
        mark options= {scale=1,draw=black,fill=mapped color},
        line join=round,
    },
    errormarks/.style={
        cycle list/Set2,
        cycle multiindex* list={
            mark list*\nextlist
            Set2\nextlist
            linestyles\nextlist
            thick
        },
        mark options={scale=0.7,line width=0.2pt},
        line join=round,
    }
}

\pgfplotsset{
    base style/.style={
        line join=round,
        unbounded coords = jump,
        width=82mm,
        height=66mm,
        grid style={line width=.2pt},
        tick style={draw=none},
        legend cell align = left,
        label style={font=\scriptsize},
        tick label style={font=\tiny},
    },
    stdplot/.style={
        base style,
        enlarge x limits = false,
        grid=both,
        stdmarks bw,
        xlabel = {Time domain},
        legend style={nodes={scale=0.8, transform shape}},
    },
    logplot/.style={
        base style,
        enlarge x limits = false,
        grid=both,
        errormarks bw,
        xlabel = {Domain},
        legend style={nodes={scale=0.5, transform shape}},
    },
    method comparison/.style={
        base style,
        enlarge x limits = false,
        grid=both,
        cycle list/Greys,
        cycle multiindex* list={
            scatter markers\nextlist
            Greys-A,Greys-C,Greys-J\nextlist
            draw=black\nextlist
        },
        mark options={scale=1.05},
    },
    method comparison four/.style={
        base style,
        enlarge x limits = false,
        grid=both,
        cycle list/Greys,
        cycle multiindex* list={
            scatter markers alt\nextlist
            Greys-A,Greys-C,Greys-E,Greys-J\nextlist
            draw=black\nextlist
        },
        mark options={scale=1.05},
    },
}

\DeclarePairedDelimiterX{\closedinterval}[2]{[}{]}{#1,#2}
\DeclarePairedDelimiterX{\norm}[1]{\lVert}{\rVert}{#1}

\newcommand{\setS}{\mathcal{S}_N}
\newcommand{\setC}{\mathcal{C}_N}
\newcommand{\xid}{\sigma}
\newcommand{\dualD}{\bm{D}^{\bm{\star}}}

\usepackage{geometry}
\title{Costate Convergence with Legendre-Lobatto Collocation for Trajectory Optimization}

\author{Jos\'e Garrido\thanks{German Aerospace Center, Institute of Space Systems, GNC Department, Robert-Hooke-Str. 7, 28359 Bremen} \\ \href{mailto:jose.valeriogarrido@dlr.de}{jose.valeriogarrido@dlr.de}
    \and Artemi Makarow\footnotemark[1] \\ \href{mailto:artemi.makarow@dlr.de}{artemi.makarow@dlr.de}
    \and Marco Sagliano\footnotemark[1] \\ \href{mailto:marco.sagliano@dlr.de}{marco.sagliano@dlr.de}
    \and David Seelbinder\footnotemark[1] \\ \href{mailto:david.seelbinder@dlr.de}{david.seelbinder@dlr.de}
\and Stephan Theil\footnotemark[1] \\ \href{mailto:stephan.theil@dlr.de}{stephan.theil@dlr.de} }

\date{}

\usepackage{amsthm}

\theoremstyle{definition}
\newtheorem{thm}{Theorem}
\theoremstyle{definition}
\newtheorem{cor}{Corollary}
\theoremstyle{definition}
\newtheorem{lem}{Lemma}
\theoremstyle{definition}
\newtheorem{prop}{Proposition}
\theoremstyle{definition}
\newtheorem{defn}{Definition}
\theoremstyle{definition}
\newtheorem{rem}{Remark}
\theoremstyle{definition}
\newtheorem{prob}{Problem}

\numberwithin{equation}{section}

\usepackage{hyperref}

\begin{document}

\maketitle

\begin{abstract}
    This paper introduces a new method of discretization that collocates both endpoints of the domain and enables the complete convergence of the costate variables associated with the Hamilton boundary-value problem.
    This is achieved through the inclusion of an \emph{exceptional sample} to the roots of the Legendre-Lobatto polynomial, thus promoting the associated differentiation matrix to be full-rank.
    We study the location of the new sample such that the differentiation matrix is the most robust to perturbations and we prove that this location is also the choice that mitigates the Runge phenomenon associated with polynomial interpolation.
    Two benchmark problems are successfully implemented in support of our theoretical findings.
    The new method is observed to converge exponentially with the number of discretization points used.
    \par\vskip\baselineskip\noindent
\end{abstract}

\section{Introduction}
Gauss-Lagrange methods, also known as pseudospectral methods, are a group of numerical methods that have seen great implementation rates over the past decades, specially applied to trajectory optimization problems which are often analytically intractable \cite{Betts1998,Rao2009,Conway2011,Rao2014a,Malyuta2021}. 

These methods are a subset of \emph{direct methods}, which enable the explicit enforcement of differential and algebraic constraints without the need to derive the costate equation associated with the Hamilton boundary-value problem.
As a result, they are practical and suitable for general purpose implementations.
In addition, the properties of the Gaussian quadrature associated with these methods provide a good trade-off between accuracy and a sparse discretization of the domain, as these methods exhibit \emph{exponential convergence rates}, see, \emph{e.g.}, \cite{Garg2011}.

The usual approach taken with these methods for trajectory optimization is to transcribe the optimal control problem into a nonlinear programming problem (NLP) and use an NLP solver to obtain a vectorized solution.
Having obtained such solution, it is then trivial to extract the individual optimal control variables.
In this context, three of the most popular Gauss-Lagrange methods are the Legendre-Lobatto method, the Legendre-Gauss method, and the Legendre-Radau method.

The Legendre-Lobatto method developed in \cite{Elnagar1995,Fahroo2001,Ross2001,Ross2003,Ross2006} is attractive because it enables the collocation of the endpoints of the domain. However, the solutions obtained with this method typically exhibit lack of convergence of the costates, see, \emph{e.g.}, \cite{Garg2011}, and as a result the accuracy of the primal variables is reduced even in a converged solution.
It was shown in \cite{Garg2011} that the associated differentiation matrix is rank-deficient and this is the leading cause for the systematic failure of this method.

In \cite{Benson2004,Rao2010} the authors propose the Legendre-Gauss method, which guarantees the convergence of the costate variables at the expense of not collocating the endpoints of the domain. In addition, the method also requires the enforcement of one quadrature constraint per state variable in any given problem, which complicates implementation.

The Legendre-Radau method introduced in \cite{Garg2011} preserves convergence of the costate variables, however, only one of the domain endpoints is collocated, therefore, control and costate variables are not obtained at the opposing endpoint.

Alternatively, in \cite{Liu2014a} the authors employ a Hermite interpolation scheme which enables a well conditioned problem with costate convergence guarantees.
However, the method still lacks the complete collocation of the terminal endpoint. As a result, the costate variables at this point are not obtained directly.

In this paper, we develop a method that enables the collocation of both endpoints of the domain and that guarantees convergence of the costate variables at every point.

By using the nodes of the Legendre-Lobatto method and including an additional discretization point at which differential constraints are not enforced, hereinafter referred to as \emph{exceptional sample}, we are able to systematically generate an appropriate differentiation matrix.
This is achieved by adhering to a reasonable qualification metric for the location of the exceptional sample and subsequently deriving a global minimizer.

As proof of concept, the new method is implemented for two classic optimal control problems for which the standard Legendre-Lobatto method has been shown to fail in the past. The precision of the new method is compared with that of the Legendre-Gauss method and of the Legendre-Radau method with a benchmark problem for which there exists an analytic solution.

\section{Preliminaries}
To obtain an adequate discretization scheme, this paper hinges on the concept of a \emph{differentiation matrix}. This motivates the following definition.
\begin{defn} \label{defn:differentiation-matrix}
    Let $\mathcal{A} = \{\tau_i \in \mathbb{R}: i=1,2,\dots,M \}$ and $\mathcal{B} = \{\tau_j \in \mathbb{R}: j=1,2,\dots,N \}$ be sets of abscissas such that $\mathcal{B} \subseteq \mathcal{A}$, then a matrix $\bm{D} \in \mathbb{R}^{N \times M}$ is called a \emph{first derivative differentiation matrix} over the set $\mathcal{A}$ if it satisfies
    \begin{equation} \label{eq:diff-defn}
        \bm{D} \bm{V} = \bm{V'}
    \end{equation}
where, for some integer $R$ such that $1 \leq R \leq M-1$, $\bm{V} \in \mathbb{R}^{M \times (R+1)}$ denotes a Vandermonde matrix over the abscissas of set $\mathcal{A}$, see, \emph{e.g.}, \cite{Trefethen1997}, and $\bm{V'} \in \mathbb{R}^{N \times (R+1)}$ denotes the respective matrix of first derivatives over the abscissas of set $\mathcal{B}$. 
    If \eqref{eq:diff-defn} holds for an integer $R$ but it does not hold for $R+1$, then $R$ denotes \emph{the order of accuracy of the method}.
    Matrices $\bm{V}$ and $\bm{V'}$ are defined, respectively, as
    \begin{align}
        \bm{V} &= \begin{bsmallmatrix}
            1 & \tau_1 & \tau_1^2 & \displaystyle \cdots & \tau_1^R \\
            1 & \tau_2 & \tau_2^2 & \displaystyle \cdots & \tau_2^R \\
            \vdots & \vdots & \vdots & \ddots & \vphantom{\int\limits^x}\vdots \\
            1 & \tau_M & \tau_M^2 & \displaystyle \cdots & \tau_M^R
        \end{bsmallmatrix}, \quad \tau_i \in \mathcal{A} \\[2ex]
        \bm{V'} &=  \begin{bsmallmatrix}
            0 & 1 & 2\tau_1  & \displaystyle \cdots & R\tau_1^{R-1} \\
            0 & 1 & 2\tau_2  & \displaystyle \cdots & R\tau_2^{R-1} \\
            \vdots & \vdots & \vdots & \ddots & \vphantom{\int\limits^x}\vdots \\
            0 & 1 & 2\tau_N  & \displaystyle \cdots & R\tau_N^{R-1}
        \end{bsmallmatrix}, \quad \tau_j \in \mathcal{B}
    \end{align}
\end{defn}

\begin{lem} \label{lem:diff-have-null-space}
    Let $\bm{D} \in \mathbb{R}^{N\times M}$ be some differentiation matrix according to Definition \ref{defn:differentiation-matrix}, then the vector $\bm{1} \in \{1\}^{M \times 1}$ is a member of the null-space of $\bm{D}$.
\end{lem}
\begin{proof}[Proof of Lemma \ref{lem:diff-have-null-space}]
    From Definition \ref{defn:differentiation-matrix} we have, necessarily, $\bm{D} \bm{1} = \bm{0}$, with $\bm{1} \in \{1\}^{M\times 1}$ and $\bm{0} \in \{0\}^{N\times 1}$.
    This implies that there exists at least one vector $\bm{z} \in \mathbb{R}^{M\times 1}$ such that $\bm{D} \bm{z} = \bm{0}$ with $\bm{z} \neq \bm{0}$.
\end{proof}

\begin{lem} \label{lem:rank-bounds}
    The rank of a differentiation matrix $\bm{D} \in \mathbb{R}^{N\times M}$ is bounded as follows:
    \begin{equation}
        \min(N, R) \leq \texttt{\upshape Rank}(\bm{D}) \leq \min( N, M-1 )\,,
    \end{equation}
    where $R$ denotes the order of accuracy of the particular method according to Definition \ref{defn:differentiation-matrix}.
\end{lem}
\begin{proof}[Proof of Lemma \ref{lem:rank-bounds}]
    From Definition \ref{defn:differentiation-matrix}, a lower bound is established from the rank of matrix $\bm{V'}$:
    \begin{equation}
        \texttt{\upshape Rank}(\bm{V'}) = \min(N, R) \leq \texttt{\upshape Rank}(\bm{D}) \,,
    \end{equation}
    where $R$ denotes the order of accuracy of the particular method of differentiation.
    An upper bound can be established from the dimensions of matrix $\bm{D}$ itself, as $\texttt{\upshape Rank}(\bm{D}) \leq \min( N, M )$.
    This boundary is narrowed further by subtracting to the number of columns the value $1$, which is the minimum dimension of the null-space of $\bm{D}$, see Lemma \ref{lem:diff-have-null-space} and \cite[Theorem 2.3]{Friedberg2003}.
\end{proof}

The class of matrices specified in Definition \ref{defn:differentiation-matrix} allows the computation of total derivates of some continuous quantity given $M$ observations of the quantity itself. We now consider the inverse problematic: obtain the sequence of values from observations of the derivates.
\begin{prob} \label{prob:discrete-initial-value-problem}
    Let $y: \mathbb{R} \mapsto \mathbb{R}$ be a continuous function, furthermore, let $\bm{y} \in \mathbb{R}^{M\times 1}$ be a discrete sequence of values $y_i$ associated with a sequence of samples of $\tau_i \in \mathbb{R}$ with $i=1,2,\dots,M$, as follows:
    \begin{equation}
        \bm{y} = \begin{bmatrix}
            y_1 & y_2 & \dots & y_M
        \end{bmatrix}^\intercal = \begin{bmatrix}
            y(\tau_1) & y(\tau_2) & \dots & y(\tau_M)
        \end{bmatrix}^\intercal \,,
    \end{equation}
    where $^\intercal$ denotes the \emph{transpose} operator.
    Determine the sequence $\bm{y}$ such that
\begin{equation} \label{eq:system-derivative-observations}
   \begin{bmatrix}
       \bm{D} \\ \bm{c}^\intercal
   \end{bmatrix}
   \bm{y} = \begin{bmatrix}
       \bm{d} \\ b
   \end{bmatrix} \,,
\end{equation}
where $\bm{D} \in \mathbb{R}^{N \times M}$ is a differentiation matrix according to Definition \ref{defn:differentiation-matrix}, $\bm{c} \in \mathbb{R}^{M \times 1}$ is a coefficient vector such that $\bm{c}^\intercal \bm{1} \neq 0$, \emph{i.e.}, $\bm{c}$ spans the null-space of $\bm{D}$, $\bm{d} \in \mathbb{R}^{N\times 1}$ is a sequence of observations of the derivative of $y(\tau)$ with respect to $\tau$ at $N$ locations, and $b \in \mathbb{R}$ is an arbitrary bias specification. 
\end{prob}

\begin{rem} \label{lem:diff-requirement}
    The linear system of Problem \ref{prob:discrete-initial-value-problem} has a unique solution if and only if matrix $\bm{D}$ consists of $N$ rows and $N+1$ columns, and if this matrix is full-rank.
\end{rem}

Remark \ref{lem:diff-requirement} motivates the following definition.

\begin{defn} \label{defn:inversion-ready}
    Let $\bm{D} \in \mathbb{R}^{N \times (N+1)}$ be a differentiation matrix according to Definition \ref{defn:differentiation-matrix} over some set of abscissas $\{\tau_i \in \mathbb{R} : i = 1,2,\ldots, N+1\}$ and for some order of accuracy $R$.
    Matrix $\bm{D}$ is an \emph{inversion-ready differentiation matrix} if $\texttt{Rank}(\bm{D}) = N$.
\end{defn}

Definition \ref{defn:inversion-ready} establishes a class of matrices which can be used to solve Problem \ref{prob:discrete-initial-value-problem}.

\section{New method of discretization}
In this section, the Lobatto polynomials are introduced and some relevant associated properties are established.
Furthermore, an approach to achieve an inversion-ready differentiation matrix according to Definition \ref{defn:inversion-ready} is proposed in this section.
At last, we assess the potential Runge's phenomenon associated with the resulting polynomial interpolation, and we briefly overview the Gaussian quadrature rule for the calculation of definite integrals with the new method.

\subsection{Lobatto polynomials}
\label{sec:polynomials}

For an integer $N \geq 0$, let $\mathcal{P}_N(\tau)$ denote the Legendre polynomial of degree $N$, with $\tau \in \closedinterval{-1}{1}$. In addition, for $N \geq 2$, let $\mathcal{L}_N(\tau)$ denote the Lobatto polynomial of degree $N$. Together, these polynomials satisfy the following equations, see \cite{Radau1880,Abramovitz1972,Szegoe1975}:
\begin{align}
    \mathcal{L}_N(\tau) &= (\tau^2-1) \dot{\mathcal{P}}_{N-1}(\tau) \label{eq:lobatto-polynomial} \,, \\
    \dot{\mathcal{L}}_N (\tau) &= N (N-1) \mathcal{P}_{N-1} (\tau) \label{eq:lobatto-first-order} \,, \\
    (\tau^2-1) \ddot{\mathcal{L}}_N (\tau) &= N (N-1) \mathcal{L}_N (\tau) \label{eq:lobatto-second-order} \,, 
\end{align}
with $ \dot{\mathcal{L}}_N (\tau) = \diff{\mathcal{L}_N(\tau)}{\tau}$ and $ \ddot{\mathcal{L}}_N (\tau) = \diff[2]{\mathcal{L}_N(\tau)}{\tau}$.

In these expressions, the integer $N$ always matches the number of roots of the respective functions, see, \emph{e.g.}, \cite{Radau1880,Hildebrand1987}.
The roots of Legendre polynomials $\mathcal{P}_N(\cdot)$ can be obtained via the methods developed in \cite{Golub1969,Bogaert2012,Bogaert2014,Glaser2007,Hale2013}, and the roots of Lobatto polynomials $\mathcal{L}_N(\cdot)$ can be obtained via the methods developed in \cite{Gautschi2000,Golub1973}.

\begin{rem} \label{rem:lobatto-zero-endpoints}
    Due to \eqref{eq:lobatto-polynomial}, the points $-1$ and $1$ are roots of $\mathcal{L}_N(\tau)$, \emph{i.e.}, $\mathcal{L}_N(1) = \mathcal{L}_N(-1) = 0$.
\end{rem}

\begin{rem} \label{rem:lobatto-extrema}
    Due to \eqref{eq:lobatto-first-order}, the stationary points of $\mathcal{L}_N(\tau)$, \emph{i.e.}, the points where $\dot{\mathcal{L}}_N(\tau) = 0$, are given by the roots of $\mathcal{P}_{N-1}(\tau)$.
\end{rem}

\begin{lem} \label{lem:lobatto-symmetry}
    Lobatto polynomials exhibit the following symmetry:
    \begin{equation} \label{eq:lobatto-symmetry}
        \mathcal{L}_N (\tau) = (-1)^N \mathcal{L}_N(-\tau) .
    \end{equation}
\end{lem}
\begin{proof}[Proof of Lemma \ref{lem:lobatto-symmetry}]
    From \cite[Eq. (4.7.4)]{Szegoe1975} the following holds:
    \begin{equation}
        \mathcal{P}_{N-1} (\tau) = (-1)^{N-1} \mathcal{P}_{N-1}(-\tau) \,.
    \end{equation}
    Taking the first derivative with respect to $\tau$, and multiplying both sides by $\tau^2 -1$ yields
    \begin{equation}
        (\tau^2-1) \dot{\mathcal{P}}_{N-1} (\tau) = (-1)^N (\tau^2-1) \dot{\mathcal{P}}_{N-1}(-\tau) \,,
    \end{equation}
    which, recalling \eqref{eq:lobatto-polynomial}, is identical to \eqref{eq:lobatto-symmetry} for $\tau \neq \pm 1$.
    The relationship still holds for $\tau = \pm 1$ due to \eqref{eq:lobatto-polynomial}, see Remark~\ref{rem:lobatto-zero-endpoints}.
\end{proof}

\subsection{Achieving an inversion-ready Lobatto differentiation matrix}
Given the set of roots of the Lobatto polynomial $\mathcal{L}_N(\cdot)$ defined in \eqref{eq:lobatto-polynomial}, the new method suggests the inclusion of \emph{an additional point of discretization along the domain at which the derivative of a continuous function of interest is not observable}, see Problem \ref{prob:discrete-initial-value-problem}.
In this section we show that this approach yields a differentiation matrix which is indeed inversion-ready according to Definition \ref{defn:inversion-ready}.

To simplify the notation henceforth, we define the set of roots of the Lobatto polynomial of degree $N$ as
\begin{equation} \label{eq:set-lobatto-roots}
    \mathcal{T}_N = \{ \tau_k \in \closedinterval{-1}{1}: \mathcal{L}_N(\tau_k) = 0 \}\,, \quad k = 1,2,\ldots,N  \,.
\end{equation}
Furthermore, let $\setS$ be a set that indexes the $N+1$ abscissas of the new method, as follows
\begin{equation} \label{eq:set-s}
    \setS = \{ 1,2, \dots, N, N+1\} \,.
\end{equation}
Then, let $\xid = N + 1$ be the \emph{index associated with the exceptional sample}, and let
\begin{equation} \label{eq:set-c}
    \setC = \setS \backslash \{\xid\} = \{ 1,2, \dots, N \}
\end{equation}
be a set that indexes the set of abscissas $\mathcal{T}_N$.

Consider a function $y: \closedinterval{-1}{1} \mapsto \mathbb{R}$.
Given $N+1$ observations of $y(\cdot)$ at $N+1$ samples $\tau_i \in \mathcal{T}_N \cup \{\tau_\xid\}$ for $i \in \setS$, where $\tau_\xid$ denotes the abscissa of the exceptional sample, an approximation of $y(\cdot)$, denoted $\tilde{y}(\cdot)$, can be constructed based on a polynomial of degree $N$ through the use of Lagrange polynomial interpolation, see, \emph{e.g.}, \cite{Berrut2004,Ross2003,Fahroo2008,Rao2010,Garg2011}, such that
\begin{equation} \label{eq:sequence-approximation}
    y(\tau) \approx \tilde{y}(\tau) = \sum_{i \in \setS} l_i(\tau) y(\tau_i) \,, \quad \tau \in \closedinterval{-1}{1} \,,
\end{equation}
where $l_i(\tau)$ denotes the Lagrange interpolating polynomial of index $i$. These polynomials are constructed as:
\begin{equation} \label{eq:poly-interp}
    l_i(\tau) = \prod_{j \in \setS \backslash \{i\}} \frac{\tau - \tau_j}{\tau_i - \tau_j} \,.
\end{equation}
From \eqref{eq:poly-interp} it is straightforward to derive that these polynomials satisfy the following:
\begin{equation} \label{eq:lagrange-kronecker-delta}
    l_i(\tau_k) = \delta_{ki} = \begin{cases} 1 \text{ if } k=i \,, \\ 0 \text{ if } k \neq i \,, \end{cases} \quad i,k \in \setS \,,
\end{equation}
where
$\tau_k \in \mathcal{T}_N \cup \{\tau_\xid\}$ denote the discretization points and $\delta_{ki}$ is called the Kronecker delta.

Finally, based on this interpolation approach, a matrix $\bm{D} \in \mathbb{R}^{N\times (N+1)}$ can be constructed as \cite{Berrut2004,Rao2010,Garg2011}
\begin{multline} \label{eq:diff-matrix}
        \bm{D}_{ki} = \dot{l}_i( \tau_k ) = \sum_{m \in \setS} \dfrac{ \prod_{j \in \setS \backslash \{i,m\}} ( \tau_k - \tau_j ) }{ \prod_{j \in \setS \backslash \{i\}} ( \tau_i - \tau_j ) } \,, \\[1ex] 
                                             k \in \setC \,, \: i \in \setS \,,
\end{multline}
where
$\tau_i \in \mathcal{T}_N \cup \{\tau_\xid\}$, 
and $\bm{D}_{ki}$ denotes the element of matrix $\bm{D}$ in row $k$ and column $i$.

\begin{thm} \label{thm:new-matrix-is-inversion-ready}
    Matrix $\bm{D} \in \mathbb{R}^{N \times (N+1)}$ specified in \eqref{eq:diff-matrix} is an inversion-ready differentiation matrix over the set of abscissas $\mathcal{T}_N \cup \{ \tau_\xid\}$ with order of accuracy $N$, according to Definition \ref{defn:inversion-ready}.
\end{thm}
\begin{proof}[Proof of Theorem \ref{thm:new-matrix-is-inversion-ready}]
    It is established that matrix $\bm{D}$ has $N$ rows and $N+1$ columns, it remains to be shown that this matrix is indeed a differentiation matrix according to Definition \ref{defn:differentiation-matrix} and that it is a full-rank matrix.
    From \cite{Davis1975,Epperson1987}, given $N+1$ distinct abscissas $\tau_i$ with $i \in \setS$, it is known that the interpolation formula expressed in \eqref{eq:sequence-approximation} is exact if the original function is a polynomial of degree $N$ or less, \emph{i.e.}, the equality $\tau^R = \sum_{i \in \setS} l_i(\tau) \tau_i^R $ holds for an integer $R \leq N$, where $l_i(\tau)$ denotes the Lagrange interpolating polynomial of index $i$. Applying the derivative with respect to $\tau$ to both sides, we get
    \begin{equation}
        \diff{\tau^R}{\tau} = R \tau^{R-1} = \sum_{i \in \setS} \dot{l}_i(\tau) \tau_i^R \,.
    \end{equation}
    Letting $\tau = \tau_k$ and recalling \eqref{eq:diff-matrix} we obtain
    \begin{equation}
        \sum_{i \in \setS} \dot{l}_i(\tau_k) \tau_i^R = \sum_{i \in \setS} \bm{D}_{ki} \tau_i^R = R \tau_k^{R-1}  \,,
    \end{equation}
    which matches Definition \ref{defn:differentiation-matrix} with order of accuracy $N$, therefore, from Lemma \ref{lem:rank-bounds}, we have $N \leq \texttt{Rank}(\bm{D}) \leq N$, implying $\texttt{Rank}(\bm{D}) = N$.
\end{proof}

\begin{rem} \label{rem:prop-is-true-for-any-abscissa}
    Theorem \ref{thm:new-matrix-is-inversion-ready} holds for all choices of $\tau_\xid$, as long as the abscissas are unique.
\end{rem}

\subsection{Determining an appropriate location for the exceptional sample}
We limit the search domain for the exceptional sample to the closed interval $\closedinterval{-1}{1}$, such that the domain of interpolation is preserved with respect to that of the unmodified roots of the Lobatto polynomial.

To determine an appropriate location $\tau_\xid$ the following metric is adopted:
\begin{equation} \label{eq:simplified-problem}
    \underset{\tau_\xid}{\text{minimize}} \quad \omega(\tau_\xid) = \norm{ \bm{D}(\tau_\xid) \bm{\delta}_\xid } \,,
\end{equation}
where
$\bm{D}(\tau_\xid)$ is an explicit parametrization of matrix $\bm{D}$ with respect to $\tau_\xid$, characterized by $\bm{D} : \closedinterval{-1}{1} \mapsto \mathbb{R}^{N\times (N+1)}$.
Furthermore, $\bm{\delta}_\xid$ is a vector of zeros except at index $\xid$ where it takes the value $\varepsilon \neq 0$, and the vector norm operator $\norm{\cdot}$ can be chosen arbitrarily for reasons that will become apparent shortly.
This metric is motivated by our intention of minimizing the degradation of interpolation quality with respect to the original set of Lobatto nodes given a perturbation at the exceptional sample.

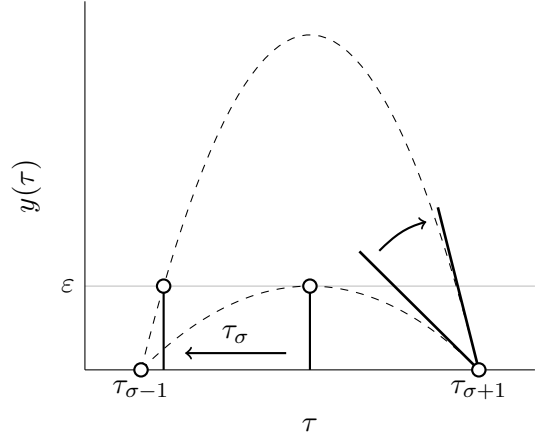
\begin{figure}[bt!]
    \centering
    \begin{tikzpicture}[trim axis left, trim axis right]
        \begin{axis}[
            width=75mm,
            ymin=0, ymax=2.2,
            ymajorgrids=true,
            axis equal, axis lines*=left,
            xtick={-1, 1},
            xticklabels={$\tau_{\xid-1}$,$\tau_{\xid+1}$},
            ytick={0.5},
            yticklabels={$\varepsilon$},
            xlabel=$\tau$,
            ylabel=$y(\tau)$,
            tick style={draw=none},
            scatter/classes={
                a={mark=*,white, draw=black,mark size = 2.5pt},
                good={mark=*,black, mark size = 2pt},
                bad={mark=*,black, mark size = 2pt}
            }
            ]
            \addplot[domain=-1:1, samples=50, dashed]{0.5-0.5*x^2};     
            \addplot[domain=-1:1, samples=50, dashed]{2-2*x^2};         
            \addplot[domain=(1 - 1/ ( (1^2+1))^(1/2)):1, samples=3,line width = 1pt]{-1*(x-1)};
            \addplot[domain=(1 - 1/ ( (4^2+1))^(1/2)):1, samples=2,line width = 1pt]{-4*(x-1)};
            \addplot[scatter, only marks, scatter src=explicit symbolic, ycomb, thick]
            coordinates {
                (-1,0) [a]
                (1,0) [a]
                (0,0.5) [a]
                (-(3/4)^(1/2),0.5) [a]
            };

            \node (source) at (axis cs:-0.08,0.10){};
            \node (destination) at (axis cs:-0.80,0.10){};
            \draw[->, thick](source)--(destination) node [midway,above] {$\tau_\xid$} ;

            \node (source) at (axis cs:0.35,0.65){};
            \node (destination) at (axis cs:0.75,0.90){};
            \draw[->,thick](source) to [out=45,in=195] (destination);
        \end{axis}
    \end{tikzpicture}
    \caption{Sensitivity of the derivative at interpolation nodes as a consequence of the location of a perturbation $\varepsilon$ at the node of index $\xid$. Poor choice of $\tau_\xid$ leads to steeper derivatives.}
    \label{fig:intuition}
\end{figure}

Figure \ref{fig:intuition} shows an example polynomial interpolation, denoted $y(\tau)$, about three points for two distinct choices of the location $\tau_\xid$ of a perturbation. The perturbation, $\varepsilon$, is kept constant for the two choices. The figure illustrates that the magnitude of the derivative at interpolation nodes is sensitive to the location of the node where the perturbation occurs.

The metric expressed in \eqref{eq:simplified-problem} is reduced to the problem of finding the minimizer of the norm of the column of matrix $\bm{D}$ with index $\xid$. In this context, from \eqref{eq:diff-matrix}, the elements in column $\xid$ of matrix $\bm{D}$ are
\begin{equation} \label{eq:D-column}
        \bm{D}_{k\xid} (\tau_\xid) =  \dfrac{\sum_{m \in \setC} \prod_{j \in \setC \backslash\{m\}} ( \tau_k - \tau_j )}{ \prod_{j \in \setC} ( \tau_\xid - \tau_j ) }  \,, \quad k \in \setC \,,
\end{equation}
where the index set $\setC$ is defined in \eqref{eq:set-c}.

It becomes apparent that the denominator and the numerator terms are proportional to the Lobatto polynomial and its derivative, respectively. Thus, \eqref{eq:D-column} can be rewritten as 
\begin{equation} \label{eq:Dke-with-lobatto}
    \bm{D}_{k\xid} (\tau_\xid) = \dfrac{\dot{\mathcal{L}}_N(\tau_k)}{ \mathcal{L}_N(\tau_\xid) } \,, \quad k \in \setC \,.
\end{equation}

Notice that the denominator term, which is the only term that depends on $\tau_\xid$, is identical for the entire column, \emph{i.e.}, it does not change with $k$. Therefore, let a constant $\alpha \in \mathbb{R}$ be such that
\begin{equation} \label{eq:alpha}
    \alpha = \norm{ \mathcal{L}_N(\tau_\xid) \bm{D}(\tau_\xid) \bm{\delta}_\xid } \,.
\end{equation}
Notice that the right-hand-side of \eqref{eq:alpha} is constant with respect to $\tau_\xid$, see \eqref{eq:Dke-with-lobatto}.
In this way, we can rewrite the objective from \eqref{eq:simplified-problem} as
\begin{equation} \label{eq:omega-final}
    \omega(\tau_\xid) = \dfrac{ \alpha }{| \mathcal{L}_N(\tau_\xid) |} \,.
\end{equation}

Finally, notice that the vector norm operator in \eqref{eq:simplified-problem} and \eqref{eq:alpha} is arbitrary. The choice of norm operator only scales the value of $\alpha$.

\begin{rem} \label{rem:maxL-minOmega}
    Due to \eqref{eq:omega-final}, the global minimizer of $\omega(\cdot)$ is given by the global maximizer of $| \mathcal{L}_N(\cdot) |$.
    Namely, let $\pi$ be such that $\omega(\pi) \leq \omega(\tau)$, $\forall \tau \in \closedinterval{-1}{1}$ holds. Then, $|\mathcal{L}_N(\pi)| \geq |\mathcal{L}_N(\tau)|$, $\forall \tau \in \closedinterval{-1}{1}$.
\end{rem}

\begin{thm} \label{thm:pattern-of-maxima}
    The global maximum of $|\mathcal{L}_N(\cdot)|$ is the stationary point of $\mathcal{L}_N(\cdot)$ whose abscissa is nearest zero.
\end{thm}
\begin{proof}[Proof of Theorem \ref{thm:pattern-of-maxima}]
    The following derivations are inspired by \cite[Proof of Theorem 7.3.1]{Szegoe1975}.
    Let us define an envelope function for $\mathcal{L}_{N}^2(\cdot)$ as
    \begin{equation} \label{eq:lobatto-envelope}
        \mathcal{F}(\tau) = \mathcal{L}_{N}^2(\tau) + \tfrac{ 1 - \tau^2 }{N(N-1)} \dot{\mathcal{L}}_{N}^2(\tau) \,,
    \end{equation}
    with $\tau \in \closedinterval{-1}{1}$,
    such that $\mathcal{F}(\tau) \geq \mathcal{L}_{N}^2(\tau)$, $\forall \tau \in \closedinterval{-1}{1}$ and $\mathcal{F}(\tau) = \mathcal{L}_{N}^2(\tau)$ when $\dot{\mathcal{L}}_{N}(\tau) = 0$ (stationary points) or when $\tau = \pm 1$ (domain endpoints).
    Applying the first derivative with respect to $\tau$ to both sides of \eqref{eq:lobatto-envelope} yields
    \begin{multline}
        \dot{\mathcal{F}}(\tau) = 2 \mathcal{L}_N(\tau) \dot{\mathcal{L}}_N (\tau) - \tfrac{2\tau}{N(N-1)} \dot{\mathcal{L}}_N^2(\tau) \\
        + \tfrac{2(1-\tau^2)}{N(N-1)}\dot{\mathcal{L}}_N (\tau)\ddot{\mathcal{L}}_N (\tau) \,,
    \end{multline}
    and factoring out the common term leads to
    \begin{multline}
        \dot{\mathcal{F}}(\tau) = 2 \dot{\mathcal{L}}_N (\tau) \bigl[ \mathcal{L}_N(\tau) - \tfrac{\tau^2-1}{N(N-1)} \ddot{\mathcal{L}}_N (\tau) \\ - \tfrac{\tau}{N(N-1)} \dot{\mathcal{L}}_N(\tau) \bigr] \,,
    \end{multline}
    and, due to \eqref{eq:lobatto-second-order}, it simplifies to
    \begin{equation}
        \dot{\mathcal{F}}(\tau) = \frac{-2\tau \dot{\mathcal{L}}_N^2(\tau)}{N(N-1)} \,.
    \end{equation}
    This shows that $\mathcal{F}(\tau)$ is increasing for $\tau<0$ and decreasing for $\tau>0$.
    Finally, let $\pi$ be such that $\dot{\mathcal{L}}_{N}(\pi) = 0$, then, from \eqref{eq:lobatto-envelope} and Lemma~\ref{lem:lobatto-symmetry} we have
    \begin{equation}
        \mathcal{F}(\pi) = \mathcal{L}_N^2 (\pi) = \mathcal{L}_N^2 (-\pi) = \mathcal{F}(-\pi) \,.
    \end{equation}
\end{proof}

\begin{cor} \label{cor:final}
    The global minimizer of $\omega(\cdot)$ is the root of $\mathcal{P}_{N-1}(\cdot)$ which is nearest zero.
\end{cor}
\begin{proof}[Poof of Corollary \ref{cor:final}]
    Due to \eqref{eq:omega-final}, see Remark \ref{rem:maxL-minOmega}, and Theorem \ref{thm:pattern-of-maxima}, the global minimizer of $\omega(\cdot)$ is the stationary point of $\mathcal{L}_N(\cdot)$ whose abscissa is closest to zero.
    Due to \eqref{eq:lobatto-first-order}, see Remark \ref{rem:lobatto-extrema}, the abscissas associated with stationary points of $\mathcal{L}_N(\cdot)$ are roots of $\mathcal{P}_{N-1}(\cdot)$.
\end{proof}

Ultimately, to minimize \eqref{eq:simplified-problem}, we select $\tau_\xid$ according to Corollary \ref{cor:final}.
This concludes the study for the location of the exceptional sample of the new method.

\subsection{Note on Runge's phenomenon}
Notably, the Lagrange interpolating polynomial of index $\xid$ is proportional to the Lobatto polynomial, in particular with a constant of proportionality of $\frac{1}{\mathcal{L}_N(\tau_\xid)}$. This result is obtained from \eqref{eq:poly-interp} as
\begin{equation}
    l_\xid(\tau) = \prod_{j \in \setC} \frac{\tau - \tau_j}{\tau_\xid - \tau_j} = \frac{\mathcal{L}_N(\tau)}{\mathcal{L}_N(\tau_\xid)} \,, \quad \tau \in \closedinterval{-1}{1} \,.
\end{equation}
Therefore, the relationship $\dot{\mathcal{L}}_N(\tau) = \dot{l}_\xid(\tau) \mathcal{L}_N(\tau_\xid)$ holds, and the stationary points of $l_\xid(\cdot)$ coincide with those of $\mathcal{L}_N(\cdot)$.
This motivates the following proposition.
\begin{prop} \label{prop:runge}
    The global minimizer of $\omega(\cdot)$ mitigates the Runge phenomenon associated with the Lagrange interpolating polynomial of index $\xid$.
\end{prop}
\begin{proof}[Proof of Proposition \ref{prop:runge}]
    Runge's phenomenon is mitigated if the interpolating polynomial remains bounded for arbitrary $N$, such that
    \begin{equation} \label{eq:runge-metric}
        |l_\xid(\tau)| \leq l_\xid(\tau_\xid) = 1 \,, \quad \forall \tau \in \closedinterval{-1}{1} \,.
    \end{equation}
    Therefore, let $\pi \in \closedinterval{-1}{1}$ be such that the inequality \eqref{eq:runge-metric} is not satisfied
    \begin{equation}
        |l_\xid(\pi)| = \frac{|\mathcal{L}_N(\pi)|}{|\mathcal{L}_N(\tau_\xid)|} > 1 \quad \Leftrightarrow \quad
        |\mathcal{L}_N(\pi)| > |\mathcal{L}_N(\tau_\xid)| \,.
    \end{equation}
    From \eqref{eq:omega-final} we know that $|\mathcal{L}_N(\tau_\xid)| \geq |\mathcal{L}_N(\tau)|\,, \quad \forall \tau \in \closedinterval{-1}{1}$, see Remark \ref{rem:maxL-minOmega}, thus:
    \begin{equation}
        |\mathcal{L}_N(\pi)| > |\mathcal{L}_N(\tau)| \,, \quad \forall \tau \in \closedinterval{-1}{1}
    \end{equation}
    which is contradictory, therefore such $\pi$ cannot exist. 
\end{proof}

\subsection{Gaussian quadrature}
Gaussian quadrature is a method of approximation of the definite integral of a function within the domain $\closedinterval{-1}{1}$ through a weighted sum of $N$ samples of the function, see, \emph{e.g.}, \cite{Radau1880,Hildebrand1987}.
In the case of Gauss-Lobatto quadrature, the method is exact for polynomial functions of degree up to $2N-3$, see, \emph{e.g.}, \cite{Radau1880,Hildebrand1987}.

In this regard, consider the polynomial function $y: \mathbb{R} \mapsto \mathbb{R}$ of degree $2N-3$ or less, and the set of roots of the Lobatto polynomial $\tau_k \in \mathcal{T}_N$ with $k \in \setC$, then the quadrature rule is expressed as follows:
\begin{equation} \label{eq:gauss-quadrature}
    \int_{-1}^{1} y(\tau) \dl \tau = \sum_{k \in \setC} w_k y(\tau_k)\,, \quad \tau_k \in \mathcal{T}_N \,,
\end{equation}
where $w_k$ are the associated Gauss-Lobatto quadrature weights which are computed as, see, \emph{e.g.}, \cite{Ross2003,Fahroo2008,Radau1880},
\begin{equation} \label{eq:weights}
    w_k = \frac{2}{N(N-1)\mathcal{P}_{N-1}^2(\tau_k) } \,, \quad k \in \setC \,.
\end{equation}

\begin{rem} \label{rem:min-degree}
    The new method employs a polynomial interpolation of degree $N$, see \eqref{eq:sequence-approximation} and \eqref{eq:poly-interp}. Therefore, to certify that the quadrature rule is always exact for this interpolation function we have to assert $2N-3 \geq N$, implying $N \geq 3$.
\end{rem}


\begin{lem} \label{lem:lagrange-inner-prod}
    Let $l_k (\tau)$ with $\tau \in \closedinterval{-1}{1}$ denote the Lagrange interpolating polynomial of index $k$, obtained with \eqref{eq:poly-interp}, furthermore, let $p(\tau)$ denote some polynomial of degree $N-3$ or less. Then the following holds:
    \begin{equation} \label{eq:lagrange-inner-product}
        \int_{-1}^{1} l_k(\tau) p(\tau) \dl \tau = w_k p (\tau_k) \,, \quad \tau_k \in \mathcal{T}_N, \: k \in \setC \,,
    \end{equation}
    where $w_k$ denotes the quadrature weight associated with index $k$, obtained from \eqref{eq:weights}.
\end{lem}
\begin{proof}[Proof of Lemma \ref{lem:lagrange-inner-prod}]
    Let $l_k(\tau)$ and $p(\tau)$ be as above. Then the following Gauss-Lobatto quadrature is exact:
    \begin{equation}
        \int_{-1}^{1} l_k(\tau) p(\tau) \dl \tau = \sum_{i \in \setC} w_i l_k(\tau_i) p(\tau_i) \,, \quad \forall k \in \setC \,,
    \end{equation}
    with $\tau_i \in \mathcal{T}_N$. Recalling \eqref{eq:lagrange-kronecker-delta} yields \eqref{eq:lagrange-inner-product}.
\end{proof}

\section{Applicability to Optimal Control}

\subsection{Problem formulation and optimality conditions}
The following is a generic optimal control problem formulation written in terms of \emph{generalized inequalities}, denoted by "$\preceq$" as in \cite{Boyd2004}. This choice is motivated by the potential applicability of the new method to convex optimal control problems.
\begin{prob} \label{prob:optimal-control}
Determine $\bm{x}(t)$, $\bm{u}(t)$, $t_0$ and $t_f$ such that
\begin{equation} \label{eq:ocp-cost}
        \Psi_0 \big( t_0, \bm{x}(t_0)\big) + \Psi_f \big( t_f, \bm{x}(t_f)\big)
                          + \int_{t_0}^{t_f} h \bigl(t, \bm{x}(t), \bm{u}(t) \bigr) \dl{t}
\end{equation}
is minimized, subject to
\begin{gather}
    \dot{\bm{x}} (t) = \bm{f} \bigl(t, \bm{x}(t), \bm{u}(t) \bigr) \label{eq:dynamics} \,, \\
    \bm{g} \bigl(t, \bm{x}(t), \bm{u}(t) \bigr)  \preceq_{\mathcal{G}} \bm{0} \,, \label{eq:path-constraints} \\
    \bm{\phi}_0 \bigl(t_0, \bm{x}(t_0) \bigr) \preceq_{\Phi_0} \bm{0} \,, \quad 
    \bm{\phi}_f \bigl(t_f, \bm{x}(t_f) \bigr) \preceq_{\Phi_f} \bm{0} \,, \label{eq:endpoint-constraints} \\
    \forall t \in \closedinterval{t_0}{t_f} \,,
\end{gather}
where $\bm{x}(t) \in \mathbb{R}^{n_x}$ and $\bm{u}(t) \in \mathbb{R}^{n_u}$ denote, respectively, the state and the control variables of the dynamic system in \eqref{eq:dynamics}.
Moreover, the constraint sets $\mathcal{G} \subseteq \mathbb{R}^{n_g}$, $\Phi_0 \subseteq \mathbb{R}^{n_{\phi,0}}$ and $\Phi_f \subseteq \mathbb{R}^{n_{\phi,f}}$ are arbitrary intersections of \emph{proper cones}, see \cite{Boyd2004}, and $\Psi_0(\cdot)$, $\Psi_f(\cdot)$, $h(\cdot)$, $\bm{f}(\cdot)$, $\bm{g}(\cdot)$, $\bm{\phi}_0(\cdot)$, and $\bm{\phi}_f(\cdot)$ are arbitrary twice-differentiable functions characterized by the following mappings:
\begin{alignat}{2}
    \Psi_0 &: \mathbb{R} \times \mathbb{R}^{n_x} &&\mapsto \mathbb{R} \,, \\
        \Psi_f &: \mathbb{R} \times \mathbb{R}^{n_x} &&\mapsto \mathbb{R} \,, \\
        h &: \mathbb{R} \times \mathbb{R}^{n_x} \times \mathbb{R}^{n_u} &&\mapsto \mathbb{R} \,, \\
        \bm{f} &: \mathbb{R} \times \mathbb{R}^{n_x} \times \mathbb{R}^{n_u} &&\mapsto \mathbb{R}^{n_x} \,, \\
        \bm{g} &: \mathbb{R} \times \mathbb{R}^{n_x} \times \mathbb{R}^{n_u} &&\mapsto \mathbb{R}^{n_g} \,, \\
        \bm{\phi}_0 &: \mathbb{R} \times \mathbb{R}^{n_x} &&\mapsto \mathbb{R}^{n_{\phi,0}} \,, \\
        \bm{\phi}_f &: \mathbb{R} \times \mathbb{R}^{n_x} &&\mapsto \mathbb{R}^{n_{\phi,f}} \,,
    \end{alignat}
where the integers $n_x$, $n_u$, $n_g$, $n_{\phi,0}$ and $n_{\phi,f}$ are non-negative and serve as arbitrary dimension specifiers.
\end{prob}

\begin{defn} \label{defn:dual-cone}
    Consider a cone $K \subseteq \mathbb{R}^{n_k}$ for some integer $n_k$. The \emph{dual cone} of $K$ is such that, see \cite{Boyd2004},
    \begin{equation}
        K^{\bm{\star}} = \{ \bm{y} \in \mathbb{R}^{n_k} : \mathbf{v}^\intercal \bm{y} \geq 0, \: \forall \mathbf{v} \in K \} \,.
    \end{equation}
\end{defn}

By introducing $\bm{\lambda}(t) \in \mathbb{R}^{n_x}$, hereinafter referred to as \emph{costate vector},
as well as 
$\bm{\mu}(t) \in \mathbb{R}^{n_g}$, 
$\bm{\nu}_0 \in \mathbb{R}^{n_{\phi,0}}$ and 
$\bm{\nu}_f \in \mathbb{R}^{n_{\phi,f}}$, 
one may define the Hamiltonian as:
\begin{multline} \label{eq:hamiltonian}
    \mathcal{H}( t, \bm{x}, \bm{u}, \bm{\lambda}, \bm{\mu}) = h ( t, \bm{x}, \bm{u} ) + \bm{\lambda}^\intercal \bm{f}( t, \bm{x}, \bm{u} ) \\
    + \bm{\mu}^\intercal \bm{g}( t, \bm{x}, \bm{u} ) \,,
\end{multline}
such that, in addition to \eqref{eq:path-constraints} and \eqref{eq:endpoint-constraints}, the solution to Problem \ref{prob:optimal-control} satisfies the following conditions
\begin{align}
    \nabla_{\bm{\lambda}} \mathcal{H}( t, \bm{x}, \bm{u}, \bm{\lambda}, \bm{\mu} ) - \dot{\bm{x}}^\intercal &= \bm{0}^\intercal \,, \label{eq:condition-lambda} \\
    \nabla_{\bm{x}} \mathcal{H}( t, \bm{x}, \bm{u}, \bm{\lambda}, \bm{\mu} ) + \dot{\bm{\lambda}}^\intercal &= \bm{0}^\intercal \,, \label{eq:condition-state} \\
    \nabla_{\bm{u}} \mathcal{H}( t, \bm{x}, \bm{u}, \bm{\lambda}, \bm{\mu} ) &= \bm{0}^\intercal \label{eq:condition-control} \,, \\
    \bm{\mu}^\intercal \bm{g}( t, \bm{x}, \bm{u} ) &= 0 \,, \label{eq:condition-slack-g} \\ 
    \bm{\nu}_0^\intercal \bm{\phi}_0 \bigl(t_0, \bm{x}(t_0) \bigr) &= 0 \,, \label{eq:condition-slack-nu-0} \\ 
    \bm{\nu}_f^\intercal \bm{\phi}_f \bigl(t_f, \bm{x}(t_f) \bigr) &= 0 \,, \label{eq:condition-slack-nu-f} \\
    \bm{\mu}(t) &\succeq_{\mathcal{G}^{\bm{\star}}} \bm{0} \\
    \bm{\nu}_0 &\succeq_{\Phi_0^{\bm{\star}}} \bm{0} \\
    \bm{\nu}_f &\succeq_{\Phi_f^{\bm{\star}}} \bm{0} \\
    \nabla_{\bm{x}} \Psi_0 \bigl(t_0, \bm{x}(t_0) \bigr) + \bm{\nu}_0^\intercal \nabla_{\bm{x}} \bm{\phi}_0 \bigl(t_0, \bm{x}(t_0) \bigr) &= -\bm{\lambda}^\intercal (t_0) \,, \label{eq:condition-state-0} \\
    \nabla_{\bm{x}} \Psi_f \bigl(t_f, \bm{x}(t_f) \bigr) + \bm{\nu}_f^\intercal \nabla_{\bm{x}} \bm{\phi}_f \bigl(t_f, \bm{x}(t_f) \bigr) &= \bm{\lambda}^\intercal (t_f) \,, \label{eq:condition-state-f} \\ 
    \diffp{\Psi_0}{t_0} \bigl(t_0, \bm{x}(t_0) \bigr) + \bm{\nu}_0^\intercal \diffp{\bm{\phi}_0}{t_0} \bigl(t_0, \bm{x}(t_0) \bigr) &= \mathcal{H}(t_0) \,,  \label{eq:condition-t0}\\
    \diffp{\Psi_f}{t_f} \bigl(t_f, \bm{x}(t_f) \bigr) + \bm{\nu}_f^\intercal \diffp{\bm{\phi}_f}{t_f} \bigl(t_f, \bm{x}(t_f) \bigr) &= - \mathcal{H}(t_f) \,, \label{eq:condition-tf}
\end{align}
where the operator $\nabla_{\mathbf{x}}$ yields the row vector of partial derivatives of the associated function with respect to the components of vector $\mathbf{x}$, and the short-hand notation $\mathcal{H}(t_0)$ stands for $\mathcal{H} \bigl( t_0, \bm{x}(t_0), \bm{u}(t_0), \bm{\lambda}(t_0), \bm{\mu}(t_0) \bigr)$. The analogous is true for $\mathcal{H}(t_f)$.
Furthermore, $\mathcal{G}^{\bm{\star}}$, $\Phi_0^{\bm{\star}}$ and $\Phi_f^{\bm{\star}}$ denote the dual cones of $\mathcal{G}$, $\Phi_0$ and $\Phi_f$, respectively, see Definition \ref{defn:dual-cone}.
In equations \eqref{eq:hamiltonian}, \eqref{eq:condition-lambda}, \eqref{eq:condition-state}, \eqref{eq:condition-control} and \eqref{eq:condition-slack-g} the parametrization of variables $\bm{x}(t)$, $\bm{u}(t)$, $\bm{\lambda}(t)$ and $\bm{\mu}(t)$ with respect to $t$ has been omitted for simplicity of notation.
Equations \eqref{eq:condition-lambda} through \eqref{eq:condition-tf} are called \emph{first order optimality conditions}.
For the formulation of the Hamiltonian and the derivation of optimality conditions see, \emph{e.g.}, \cite{Bryson1975,Kirk2004,Boyd2004}.

\subsection{Algebraic transcription and constraint enforcement}
To employ the new method we discretize Problem \ref{prob:optimal-control} according to the following strategy.
\begin{description}
    \item[State Variables]{Sampled according to $\mathcal{T}_N \cup \{\tau_\xid\}$, with $\tau_\xid$ given by Corollary \ref{cor:final}. Forming a set of $N+1$ points that are indexed with $\setS$.}
    \item[Control Variables]{Sampled according to $\mathcal{T}_N$. Forming a set of $N$ points that are indexed with $\setC$.}
\end{description}
Please recall that $\mathcal{T}_N$ is the set of roots of the Lobatto polynomial defined in \eqref{eq:set-lobatto-roots}, and that the index sets $\setS$ and $\setC$ are defined in \eqref{eq:set-s} and \eqref{eq:set-c}, respectively.

Note that by discretizing the problem we reduce the dimensionality of the search space to a finite set, therefore, there is always an intrinsic error present between discrete and analytic variables.
Regardless, for the sake of simplicity and intuitiveness, we shall use identical variable naming for respective quantities.

Given the domain of interest $\closedinterval{t_0}{t_f}$ and the normalized domain $\tau \in \closedinterval{-1}{1}$, a \emph{one-to-one} relationship can be established through the linear mapping $t : \closedinterval{-1}{1} \mapsto \closedinterval{t_0}{t_f}$, defined as:
\begin{equation} \label{eq:domain-mapping}
    t(\tau) = \frac{t_f-t_0}{2} \tau + \frac{t_f+t_0}{2} \,.
\end{equation}

In this way, we are able to enforce the path constraints specified in \eqref{eq:path-constraints} as:
\begin{equation} \label{eq:discrete-path-constraints}
    \bm{g}( t_k, \bm{x}_k, \bm{u}_k ) \preceq_{\mathcal{G}} \bm{0} \,, \quad k \in \setC \,,
\end{equation}
where $t_k = t( \tau_k )$ according to \eqref{eq:domain-mapping}, $\bm{x}_k$ and $\bm{u}_k$ denote $\bm{x}(t_k)$ and $\bm{u}(t_k)$, respectively, and $\tau_k \in \mathcal{T}_N$ with $k \in \setC$.

Moreover, recalling the Gaussian quadrature from \eqref{eq:gauss-quadrature} and the domain mapping in \eqref{eq:domain-mapping}, one may approximate the definite integral of Problem \ref{prob:optimal-control} with the following algebraic form, see, \emph{e.g.}, \cite{Garg2011,Patterson2014,Garrido2021,GarridoMScThesis}:
\begin{equation} \label{eq:algebraic-integral}
    \int_{t_0}^{t_f} h \bigl(t, \bm{x}(t), \bm{u}(t) \bigr) \dl{t} \approx \frac{t_f - t_0}{2} \sum_{k \in \setC} w_k h( t_k, \bm{x}_k, \bm{u}_k ) \,.
\end{equation}

Finally, the differential constraints expressed in \eqref{eq:dynamics} are enforced with the following algebraic construction: 
\begin{equation} \label{eq:diff-scheme}
    \bm{f}(t_k, \bm{x}_k, \bm{u}_k) - \frac{2}{t_f-t_0} \sum_{i \in \setS} \bm{D}_{ki} \bm{x}_i = \bm{0} \,, \quad k \in \setC \,,
\end{equation}
where $\bm{D}$ denotes the differentiation matrix of the new method, see Theorem \ref{thm:new-matrix-is-inversion-ready}.

\subsection{Karush–Kuhn–Tucker conditions}
In this section, let $1$ and $N$ be the indices associated with the initial and final domain endpoints, respectively.
In this way, we may define a \emph{KKT Hamiltonian} as:
\begin{multline} \label{eq:kkt-hamiltonian}
    \tilde{H}( t_0, t_f, \bm{x}_k, \bm{u}_k, \bm{\lambda}_k, \bm{\mu}_k, \bm{\nu}_0, \bm{\nu}_f ) = \\
    \tfrac{t_f-t_0}{2} w_k \big(
    h_k
    + \bm{\lambda}_k^\intercal \bm{f}_k 
    + \bm{\mu}_k^\intercal \bm{g}_k 
    \big) \\
    + \big( \Psi_0(t_0, \bm{x}_1) + \bm{\nu}_0^\intercal \bm{\phi}_0(t_0, \bm{x}_1) \big) \delta_{1k} \\
    + \big( \Psi_f(t_f, \bm{x}_N) + \bm{\nu}_f^\intercal \bm{\phi}_f(t_f, \bm{x}_N) \big) \delta_{Nk} \,,
\end{multline}
where
$\bm{\lambda}_k$ and $\bm{\mu}_k$ denote $\bm{\lambda}(t_k)$ and $\bm{\mu}(t_k)$, respectively, with $k \in \setC$.
Moreover, $\delta_{ik}$ denotes the Kronecker delta defined in \eqref{eq:lagrange-kronecker-delta},
and the short-hand $\mathbf{f}_k = \mathbf{f}\bigl( t_k, \bm{x}(t_k), \bm{u}(t_k) \bigr)$ is used for functions $h$, $\bm{f}$ and $\bm{g}$ of Problem \ref{prob:optimal-control}.

To simplify notation, we introduce the following short-hand for the KKT Hamiltonian function:
\begin{equation}
    \tilde{H}_k = \tilde{H}( t_0, t_f, \bm{x}_k, \bm{u}_k, \bm{\lambda}_k, \bm{\mu}_k, \bm{\nu}_0, \bm{\nu}_f ) \,,
\end{equation}
where the initial and final domain endpoints, $t_0$ and $t_f$, as well as the multipliers, $\bm{\nu}_0$ and $\bm{\nu}_f$, do not change for $k \in \setC$.

Using the KKT Hamiltonian, we construct an \emph{augmented cost function} as follows:
\begin{multline} \label{eq:augmented-cost}
    J( t_0, t_f, \bm{x}_j, \bm{u}_m, \bm{\lambda}_m, \bm{\mu}_m, \bm{\nu}_0, \bm{\nu}_f ) = \\
    \sum_{k \in \setC} \Bigl\{ \tilde{H}_k - w_k \bm{\lambda}_k^\intercal \sum_{i \in \setS} \bm{D}_{ki} \bm{x}_i \Bigr\}\,, \\ \quad j \in \setS, \: m \in \setC,
\end{multline}
In this way, the augmented cost in \eqref{eq:augmented-cost} is analogous to the Lagrangian function in the context of constrained parametric optimization with respect to the transcription strategy expressed in \eqref{eq:algebraic-integral}, \eqref{eq:discrete-path-constraints}, and \eqref{eq:diff-scheme}, see, \emph{e.g.}, \cite{Bertsekas2016,Boyd2004,Garg2011}.

As an intermediate step towards the KKT conditions, we highlight the gradient of the augmented cost with respect to the discrete state variables below:
\begin{equation} \label{eq:kkt-x-intermediate}
    \nabla_{\bm{x}_j} J(\cdot) = \sum_{k \in \setC} \bigl\{ \nabla_{\bm{x}_k} \tilde{H}_k \delta_{kj} - w_k \bm{\lambda}_k^\intercal \bm{D}_{kj} \bigr\} \,, \quad j \in \setS \,.
\end{equation}
Equation \eqref{eq:kkt-x-intermediate} exhibits a \emph{left} multiplication between the costate samples, $\bm{\lambda}_k$ with $k \in \setC$, and the differentiation matrix $\bm{D}$.
To convert this into a \emph{right} multiplication and obtain a linear combination in column space, we define an auxiliary matrix $\dualD \in \mathbb{R}^{N \times N}$ as:
\begin{equation} \label{eq:dual-matrix}
    \dualD_{ki} = \frac{\delta_{ki}}{w_k}( \delta_{Nk} - \delta_{1k} ) - \frac{w_i}{w_k} \bm{D}_{ik} \,, \quad k, i \in \setC\,,
\end{equation}
where $w_k$ are the quadrature weights obtained from \eqref{eq:weights}.
In this way, the KKT conditions, which are obtained by setting the gradient of $J(\cdot)$ with respect to each variable to zero, are expressed as follows:
\begin{align}
    \nabla_{\bm{\lambda}_k} \tilde{H}_k - w_k \sum_{i \in \setS} \bm{D}_{ki} \bm{x}_i^\intercal &= \bm{0}^\intercal \,, \label{eq:kkt-l} \\
    \nabla_{\bm{x}_k} \tilde{H}_k + w_k \sum_{i \in \setC} \dualD_{ki} \bm{\lambda}_i^\intercal &= \bm{\lambda}_N^\intercal \delta_{Nk} - \bm{\lambda}_1^\intercal \delta_{1k} \,, \label{eq:kkt-x} \\
    \sum_{i \in \setC} w_i \bm{\lambda}_i \bm{D}_{i\xid} &= \bm{0} \,, \label{eq:kkt-x-sigma} \\
    \nabla_{\bm{u}_k} \tilde{H}_k &= \bm{0}^\intercal \,, \label{eq:kkt-u} \\
    \bm{\mu}_k^\intercal \bm{g}_k &= 0 \,, \\ 
    \bm{\nu}_0^\intercal \bm{\phi}_0 (t_0, \bm{x}_1 ) &= 0 \,, \\ 
    \bm{\nu}_f^\intercal \bm{\phi}_f (t_f, \bm{x}_N ) &= 0 \,, \\
    \sum_{i \in \setC} \diffp{\tilde{H}_i}{t_0} = \sum_{i \in \setC} \diffp{\tilde{H}_i}{t_f} &= 0 \,,
\end{align}
with $k \in \setC$, see, \emph{e.g.}, \cite{Ross2003,Benson2006,Fahroo2008,Garg2011}.
Note that \eqref{eq:kkt-x} is a consequence of \eqref{eq:kkt-x-intermediate} for $j \in \setC$.
Similarly, \eqref{eq:kkt-x-sigma} is a consequence of \eqref{eq:kkt-x-intermediate} for $j = \xid$.

\begin{thm} \label{thm:dual-matrix-is-diff-matrix}
    Matrix $\dualD \in \mathbb{R}^{N\times N}$ specified in \eqref{eq:dual-matrix} is a differentiation matrix over the set of abscissas $\mathcal{T}_N$ with order of accuracy $N-2$, according to Definition \ref{defn:differentiation-matrix}.
\end{thm}
\begin{proof}[Proof of Theorem \ref{thm:dual-matrix-is-diff-matrix}]
    Consider a polynomial $p(\tau)$ of degree $N_p$, with $\tau \in \closedinterval{-1}{1}$ sampled at the nodes $\tau_i \in \mathcal{T}_N$, with $i \in \setC$. From \eqref{eq:dual-matrix} and \eqref{eq:diff-matrix}, we apply the matrix $\dualD$ to the samples $p(\tau_i)$, as follows:
    \begin{multline} \label{eq:dual-D-with-poly}
        \sum_{i \in \setC} \dualD_{ki} p(\tau_i) = \frac{\delta_{Nk} - \delta_{1k}}{w_k} \sum_{i \in \setC} \delta_{ki} p(\tau_i) \\ 
        - \frac{1}{w_k} \sum_{i \in \setC} w_i \dot{l}_k (\tau_i) p(\tau_i)
    \end{multline}
    where $l_k(\cdot)$ denotes the Lagrange interpolating polynomial of index $k$ defined in \eqref{eq:poly-interp}. As a result, the polynomial $\dot{l}_k(\cdot)$ has degree $N-1$. In this context, let $N_p$ be such that $(N-1) + N_p \leq 2N - 3$ leading to $N_p \leq N - 2$, then the Gaussian quadrature rule in \eqref{eq:dual-D-with-poly} is exact, allowing us to write the right-hand-side of \eqref{eq:dual-D-with-poly} as:
    \begin{equation}
         \frac{\delta_{Nk} - \delta_{1k}}{w_k} p(\tau_k) - \frac{1}{w_k} \int_{-1}^{1} \dot{l}_k (\tau) p(\tau) \dl \tau \,.
    \end{equation}
    Performing integration by parts yields
    \begin{multline}
        \frac{\delta_{Nk} - \delta_{1k}}{w_k} p(\tau_k) +
        \frac{1}{w_k} \int_{-1}^{1} l_k (\tau) \dot{p}(\tau) \dl \tau \\
         - \frac{1}{w_k}\bigl[ l_k (\tau_N) p(\tau_N) -  l_k (\tau_1) p(\tau_1)\bigr] \,.
    \end{multline}
    From \eqref{eq:lagrange-kronecker-delta}, it is apparent that the endpoint terms vanish if $k \notin \{ 1,\, N \}$, and that these terms cancel-out if $k \in \{ 1, \, N \}$.
    Furthermore, from Lemma \ref{lem:lagrange-inner-prod}, the integral term can be simplified, yielding
    \begin{equation}
            \sum_{i \in \setC} \dualD_{ki} p(\tau_i) = \frac{1}{w_k} w_k \dot{p}(\tau_k) = \dot{p}(\tau_k) \,.
    \end{equation}
\end{proof}

\begin{thm} \label{thm:costate-degree}
    The optimal sequence of costate variables, composed of $N$ total samples, $\bm{\lambda}_i$ with $i \in \setC$, is interpolatable with a polynomial of degree $N-2$ or less.
\end{thm}
\begin{proof}[Proof of Theorem \ref{thm:costate-degree}]
    Let $q: \closedinterval{-1}{1} \mapsto \mathbb{R}$ be a polynomial of degree $N-1$ for some $N\geq 2$, then we can write $q(\cdot)$ in terms of Legendre polynomials as follows:
    \begin{equation} \label{eq:q-with-Legendre}
        q(\tau) = \sum_{j = 0}^{N-1} c_j \mathcal{P}_j (\tau) = \sum_{j = 0}^{N-2} c_j \mathcal{P}_j (\tau) + c_{N-1} \mathcal{P}_{N-1} (\tau) \,,
    \end{equation}
    where $c_j$ with $j =0,\dots,N-1$ are arbitrary scalar coefficients.
    Recalling \eqref{eq:Dke-with-lobatto} and \eqref{eq:lobatto-first-order}, we can rewrite the necessary optimality condition \eqref{eq:kkt-x-sigma} in terms of $q(\tau_i)$ instead of $\bm{\lambda}_i$ as follows:
    \begin{equation} \label{eq:kkt-x-sigma-alt}
        \sum_{i \in \setC} w_i q(\tau_i)  \bm{D}_{i\xid}
        = \tfrac{N(N-1)}{\mathcal{L}_N(\tau_\xid)}
        \sum_{i \in \setC} w_i q(\tau_i) \mathcal{P}_{N-1}(\tau_i) = 0
    \end{equation}
    Eliminating the leading quotient, expanding $q(\tau_i)$ with \eqref{eq:q-with-Legendre} and replacing the resulting Gaussian quadrature rules with definite integrals when possible, yields
    \begin{multline} \label{eq:q-expanded}
        \sum_{j = 0}^{N-2} c_j \overbrace{\int_{-1}^{1} \mathcal{P}_j(\tau) \mathcal{P}_{N-1}(\tau) \dl \tau }^{= 0} \\
         + c_{N-1} \sum_{i \in \setC} w_i \mathcal{P}_{N-1}^2 (\tau_i) = 0 \,.
    \end{multline}
    In this expression, the definite integral terms vanish due to the orthogonality property of Legendre polynomials \cite{Szegoe1975}.
    It becomes apparent that \eqref{eq:q-expanded}, and consequently \eqref{eq:kkt-x-sigma-alt}, can only be satisfied if $c_{N-1} = 0$.
\end{proof}

\subsection{Covector mapping}
It has been shown that there is a one-to-one correspondence between the discrete dual variables of the optimal control problem (with respect to which the KKT conditions above have been written) and the Lagrange multipliers that are obtained from the solver in use, see, \emph{e.g.}, \cite{Ross2001,Ross2003,Fahroo2008,Garg2011,Benson2006,GarridoMScThesis}. For the augmented cost in \eqref{eq:augmented-cost} and the KKT Hamiltonian definition in \eqref{eq:kkt-hamiltonian}, this mapping is given as:
\begin{align}
    \bm{\lambda}_k &= \frac{2}{w_k(t_f - t_0)} \tilde{\bm{\lambda}}_k \,, \quad k \in \setC \\
    \bm{\mu}_k &= \frac{2}{w_k(t_f - t_0)} \tilde{\bm{\mu}}_k \,, \quad k \in \setC \\
    \bm{\nu}_0 &= \tilde{\bm{\nu}}_0 \\
    \bm{\nu}_f &= \tilde{\bm{\nu}}_f
\end{align}
where $\tilde{\bm{\lambda}}_k$, $\tilde{\bm{\mu}}_k$, $\tilde{\bm{\nu}}_0$, and $\tilde{\bm{\nu}}_f$ denote the Lagrange multipliers which are respectively associated with the differential constraints, path constraints and initial and terminal endpoint constraints.

\section{Numerical results}
In this section we implement the new method for two classical optimal control problems.
These are the low-thrust orbit raising problem \cite{Fahroo2008,Garg2011,Bryson1975} and the scalar nonlinear initial value problem \cite{Garg2011}, the later of which has a known analytic solution.
The complex step method for partial derivates is used with a step size of \num{1e-12}, see, \emph{e.g.}, \cite{Alonso2003}, and the NLP solver IPOPT \cite{Waechter2005} was employed.

Table \ref{tab:ipopt-options} shows the solver settings used for all numerical results that follow.
The tolerance values shown in Table~\ref{tab:ipopt-options} directly enforce the precision with which the KKT conditions are satisfied.
\begin{table}[bt!]
    \centering
    \noindent
    \caption{List of IPOPT options used for the benchmark problems.}
    \label{tab:ipopt-options}
    \vspace{1mm}
    \begin{tabular}{lc}
        \toprule
        Option & Value \\ \midrule
        \texttt{tol} & \num{1e-9} \\
        \texttt{dual\_inf\_tol} & \num{1e-9} \\
        \texttt{acceptable\_tol} & \num{1e-9} \\
        \texttt{max\_iter} & \num{1e3} \\
        \texttt{hessian\_approximation} & \texttt{'limited-memory'} \\
        \texttt{linear\_solver} & \texttt{'mumps'} \\
        \bottomrule
    \end{tabular}
\end{table}

\subsection{Time-invariant orbit raising problem}
The orbit raising problem \cite{Bryson1975,Fahroo2008,Garg2011} is concerned with the fixed-time transfer between two circular orbits such that the radius of the final orbit is maximized.
Here it has been adapted such that the mass of the system is included as a state, and therefore the dynamic equations are time-invariant.
The problem is stated as follows.
\begin{prob} \label{prob:orbit-raising}
    Determine $r(t)$, $\theta(t)$, $v_r(t)$, $v_\theta(t)$, $m(t)$ and $\beta(t)$, $\forall t \in \closedinterval{0}{t_f}$ such that
\begin{equation}
     -r(t_f)
\end{equation}
is minimized, subject to the equations of motion
\begin{align}
    \dot{r}(t) &= v_r(t) \,, \\ 
    \dot{\theta}(t) &= \tfrac{v_\theta(t)}{r(t)} \,, \\ 
    \dot{v}_r(t) &= \tfrac{v_\theta^2(t)}{r(t)} - \tfrac{\mu}{r^2(t)} + \tfrac{T}{m(t)} \sin \beta(t) \,,   \\
    \dot{v}_\theta(t) &= - \tfrac{v_r(t) v_\theta(t)}{r(t)} + \tfrac{T}{m(t)} \cos \beta(t) \,, \\ 
    \dot{m}(t) &=-\alpha \,, 
\end{align}
the initial conditions
\begin{align}
    r(0) &= 1 \,, \\
    \theta(0) &= 0 \,, \\
    v_r(0) &= 0 \,, \\
    v_\theta(0) &= 1 \,, \\
    m(0) &= 1 \,,
\end{align}
and the terminal conditions
\begin{align}
    v_r(t_f) &= 0 \,, \\
    v_\theta(t_f) - \sqrt{\tfrac{\mu}{r(t_f)}} &= 0 \,,
\end{align}
where $t_f = \num{3.32}$, $T = \num{0.1405}$, $\mu =\num{1}$, and $\alpha = \num{0.0749}$.
\end{prob}

Problem \ref{prob:orbit-raising} is solved with the new method for $N=\num{25}$. In this context, Fig. \ref{fig:orbit-raising-state} and Fig. \ref{fig:orbit-raising-costate} show the obtained trajectories of state and costate, respectively. The profiles are in agreement with literature, see, \emph{e.g.}, \cite{Garg2011,Fahroo2008}.
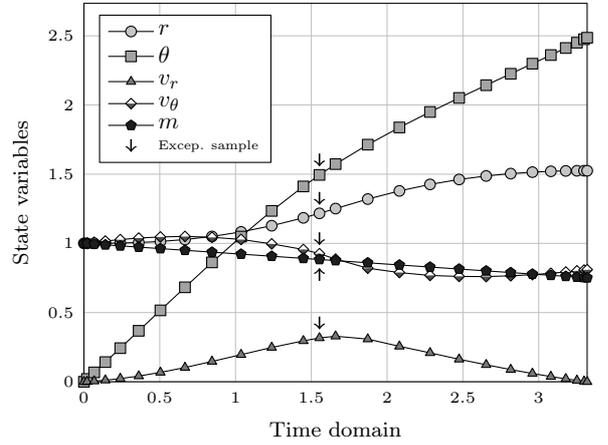
\begin{figure}[bt!]
    \raggedleft
    \begin{tikzpicture}[trim axis left, trim axis right]
    \begin{axis}[stdplot,
        ylabel = {State variables},
        ymin = 0,
        legend pos = north west
        ]

        \addplot table [x=t, y=r] {data/orbit-raising/state.txt};
        \addplot table [x=t, y=a] {data/orbit-raising/state.txt};
        \addplot table [x=t, y=vr] {data/orbit-raising/state.txt};
        \addplot table [x=t, y=va] {data/orbit-raising/state.txt};
        \addplot table [x=t, y=m] {data/orbit-raising/state.txt};

        \legend{$r$,$\theta$,$v_r$,$v_\theta$,$m$}

        \addplot[mark=arrow, only marks] coordinates {
            (1.55366555784807, 1.2169848854907+0.1)
            (1.55366555784807, 1.49474348532068+0.1)
            (1.55366555784807, 0.316600523420509+0.1)
            (1.55366555784807, 0.925598410890248+0.1)
            };
        \addlegendentry{\tiny Excep. sample}

        \addplot[mark=arrow, only marks, mark options={rotate=180}] coordinates { (1.55366555784807, 0.883630449717174-0.1) };

    \end{axis}
    \end{tikzpicture} \enspace
    \caption{State profiles for Problem \ref{prob:orbit-raising}.}
    \label{fig:orbit-raising-state}
\end{figure}
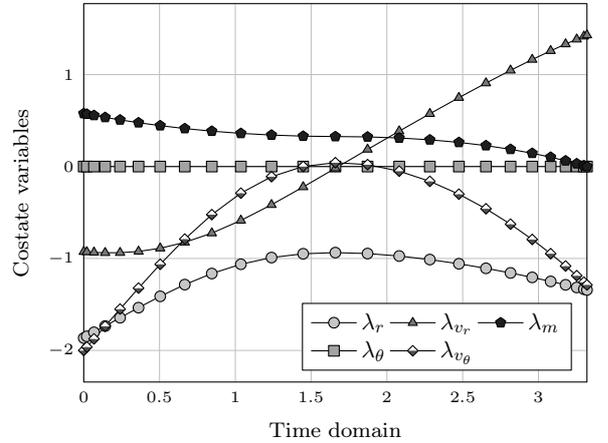
\begin{figure}[bt!]
    \raggedleft
    \begin{tikzpicture}[trim axis left, trim axis right]
    \begin{axis}[ stdplot,
        ylabel = {Costate variables},
        legend pos = south east,
        legend columns=2,
        transpose legend,
        ]
        \addplot table [x=t, y=lr] {data/orbit-raising/costate.txt};
        \addplot table [x=t, y=la] {data/orbit-raising/costate.txt};
        \addplot table [x=t, y=lvr] {data/orbit-raising/costate.txt};
        \addplot table [x=t, y=lva] {data/orbit-raising/costate.txt};
        \addplot table [x=t, y=lm] {data/orbit-raising/costate.txt};

        \legend{$\lambda_r$,$\lambda_\theta$,$\lambda_{v_r}$,$\lambda_{v_\theta}$, $\lambda_m$}
    \end{axis}
    \end{tikzpicture} \enspace
    \caption{Costate profiles for Problem \ref{prob:orbit-raising}.}
    \label{fig:orbit-raising-costate}
\end{figure}

\subsection{Convergence behaviour for an analytic problem}

In order to test the convergence rate of the new method, we shall consider a problem for which the solution is known. In this context we shall take a scalar nonlinear initial value problem from \cite{Garg2011} and \cite{Garg2010}, which is stated as follows.
\begin{prob} \label{prob:convergence}
    Determine $x(t)$ and $u(t)$, $\forall t \in \closedinterval{0}{2}$ such that
    \begin{equation}
        -x(2)
    \end{equation}
    is minimized, subject to the differential constraint
    \begin{equation}
        \dot{x}(t) = \tfrac{5}{2}( x(t) u(t) - x(t) - u^2(t) ) \,,
    \end{equation}
    and the initial condition $x(0) = 1$.
\end{prob}

The solution to Problem \ref{prob:convergence} is as follows, see \cite{Garg2011,Garg2010}:
\begin{align}
    x^*(t) &= 4/( 1 + 3 \exp{\tfrac{5t}{2}} )\\
    u^*(t) &= x^*(t)/2 \\
    \lambda^*(t) &= - \frac{\exp{\bigl(2 \ln(1 + 3 \exp{(\tfrac{5t}{2})}) - \tfrac{5t}{2}\bigr)} }{6 + 9 \exp(5) + \exp(-5) } 
\end{align}

We now show the achieved precision of the new method with respect to the known solution for many degrees of polynomial interpolation $N$ and we compare the result with the alternative Gauss-Lagrange methods, namely, the Gauss method \cite{Benson2004,Rao2010}, the Radau method \cite{Garg2011,Patterson2014} and the standard Lobatto method \cite{Ross2003,Fahroo2008}.

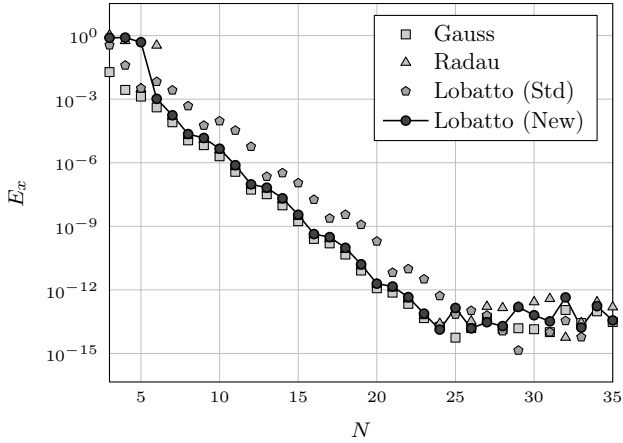
\begin{figure}[bt!]
    \raggedleft
    \begin{tikzpicture}[trim axis left, trim axis right]
    \begin{semilogyaxis}[
        method comparison four,
        only marks,
        xlabel = {$N$},
        ylabel = {$E_x$},
        legend style={nodes={scale=0.80, transform shape}},
        legend pos = north east,
        xtick distance=5,
        ytick={10^0,10^-3,10^-6,10^-9,10^-12,10^-15},
        ]
        \addplot table [x=n, y=gauss] {data/convergence/state.txt};
        \addplot table [x=n, y=radau]  {data/convergence/state.txt};
        \addplot table [x=n, y=stdlobatto] {data/convergence/state.txt};
        \addplot+[sharp plot, semithick] table [x=n, y=lobatto] {data/convergence/state.txt};

        \legend{Gauss,Radau,Lobatto (Std),Lobatto (New)}
    \end{semilogyaxis}
    \end{tikzpicture} \enspace
    \caption{State error vs. number of points used for Problem \ref{prob:convergence} with different discretization methods.}
    \label{fig:convergence-state}
\end{figure}

Figures \ref{fig:convergence-state}, \ref{fig:convergence-control}, and \ref{fig:convergence-costate} show the maximum absolute value of the errors obtained for state, control and costate, respectively, in relation to the known solution \emph{at the points where differential constraints are enforced}. These quantities are computed as follows
\begin{align}
    E_x &= \max_{k \in \setC}{\lvert \hat{x}_k - x^*(t_k) \rvert} \\
    E_u &= \max_{k \in \setC}{\lvert \hat{u}_k - u^*(t_k) \rvert} \\
    E_\lambda &= \max_{k \in \setC}{\lvert \hat{\lambda}_k - \lambda^*(t_k) \rvert} \,,
\end{align}
where $\hat{x}_k$, $\hat{u}_k$, and $\hat{\lambda}_k$ are, respectively, the discrete solutions of state, control and costate obtained with each method.
\begin{figure}[bt!]
    \raggedleft
    \begin{tikzpicture}[trim axis left, trim axis right]
    \begin{semilogyaxis}[
        method comparison four,
        only marks,
        xlabel = {$N$},
        ylabel = {$E_u$},
        legend style={nodes={scale=0.80, transform shape}},
        legend pos = north east,
        xtick distance=5,
        ]
        \addplot table [x=n, y=gauss] {data/convergence/control.txt};
        \addplot table [x=n, y=radau]  {data/convergence/control.txt};
        \addplot table [x=n, y=stdlobatto] {data/convergence/control.txt};
        \addplot+[sharp plot, semithick] table [x=n, y=lobatto] {data/convergence/control.txt};

        \legend{Gauss,Radau,Lobatto (Std),Lobatto (New)}
    \end{semilogyaxis}
    \end{tikzpicture} \enspace
    \caption{Control error vs. number of points used for Problem \ref{prob:convergence} with different discretization methods.}
    \label{fig:convergence-control}
\end{figure}
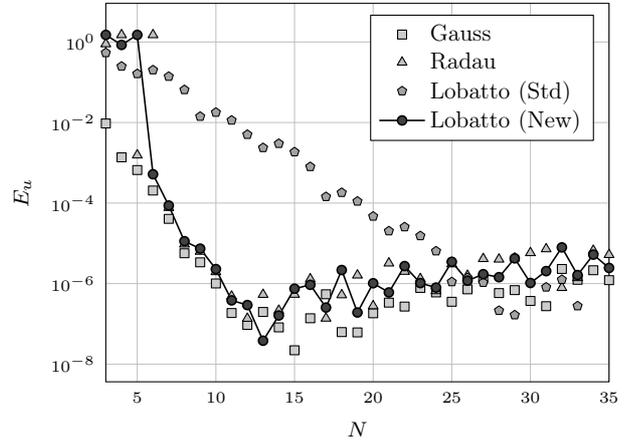
\begin{figure}[bt!]
    \raggedleft
    \begin{tikzpicture}[trim axis left, trim axis right]
    \begin{semilogyaxis}[
        method comparison four,
        only marks,
        xlabel = {$N$},
        ylabel = {$E_\lambda$},
        legend style={nodes={scale=0.80, transform shape},at={(0.97,0.6)},anchor=north east},
        ytick={10^0,10^-3,10^-6,10^-9,10^-11},
        xtick distance=5,
        ]
        \addplot table [x=n, y=gauss] {data/convergence/costate.txt};
        \addplot table [x=n, y=radau]  {data/convergence/costate.txt};
        \addplot table [x=n, y=stdlobatto] {data/convergence/costate.txt};
        \addplot+[sharp plot, semithick] table [x=n, y=lobatto] {data/convergence/costate.txt};

        \legend{Gauss,Radau,Lobatto (Std),Lobatto (New)}
    \end{semilogyaxis}
    \end{tikzpicture} \enspace
    \caption{Costate error vs. number of points used for Problem \ref{prob:convergence} with different discretization methods.}
    \label{fig:convergence-costate}
\end{figure}
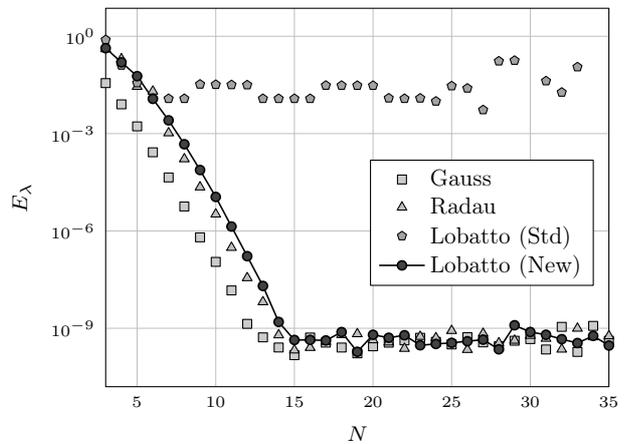

As shown in previous research, see \cite{Garg2011,Garg2010}, the standard Lobatto method, which employs a square differentiation matrix, is, in general, the most inaccurate.
In addition, note that the IPOPT solver failed to reach optimality when using the standard Lobatto method for $N=\{30,34,35\}$. The absence of the respective points is evident in Fig. \ref{fig:convergence-costate}.
In contrast, we note that the new method exhibits identical precision as the method of Radau with regards to the state and control solutions.


%
%

\section{Conclusion}
In this paper, we have developed a discretization method that enables algebraic differentiation and integration over the domain $\closedinterval{-1}{1}$.
As discretization points, the new method suggests using the roots of the Lobatto polynomial of degree $N$ as well as the root of the Legendre polynomial of degree $N-1$ that is nearest zero, the latter denoted as \emph{exceptional sample}.
We have shown that the resulting differentiation matrix exhibits important beneficial properties that classify it as \emph{inversion-ready}.
We focused on the numerical solution of optimal control problems and we saw that the new method exhibits comparable performance with respect to alternative state of the art methods.

The application of the new method to trajectory optimization is of particular interest, as it enables the direct attainment of the costate variables at both endpoints of the domain without sacrificing convergence behaviour.

\bibliographystyle{unsrt}
\bibliography{optimal-control}

\end{document}